\documentclass[11pt]{article}

\usepackage{amssymb}
\usepackage{amsmath}
\usepackage{theorem}
\usepackage{epsfig}
\usepackage{verbatim}
\usepackage{graphicx}

\newtheorem{thm}{Theorem}[section]

\newcommand{\R}{\Bbb{R}}

\newcommand{\T}{\mathbb{T}}
\newcommand{\D}{\displaystyle}

\newcommand{\grad}{\nabla}

\newcommand{\al}{\alpha}

\usepackage{color}

\newcommand{\eqdef}{\overset{\mbox{\tiny{def}}}{=}}

\begin{document}

\author{Francisco Gancedo and Robert M. Strain}

\title{Absence of splash singularities for SQG sharp fronts\\ and the Muskat problem}

\date{September 13, 2013}

\maketitle

\begin{abstract}
In this paper for either the sharp front Surface Quasi-Geostrophic equation or the Muskat problem we rule out the ``splash singularity'' blow-up scenario; in other words we prove that the contours evolving from either of these systems can not intersect at a single point while the free boundary remains smooth. Splash singularities have been shown to hold for the free boundary incompressible Euler equation in the form of the water waves contour evolution problem \cite{ADCPJ}. Our result confirms the numerical simulations in \cite{CFMR} where it is shown that the curvature blows up due to the contours collapsing at a point. Here we prove that maintaining control of the curvature will remove the possibility of pointwise interphase collapse. Another conclusion that we provide is a better understanding of the work \cite{DP3} in which squirt singularities are ruled out; here a positive volume of fluid between the contours can not be ejected in finite time.
\end{abstract}



\setcounter{tocdepth}{1}

\section{Introduction}
We consider the following general transport evolution equation:
\begin{align}
\begin{split}\label{ITE}
\zeta_t(x,t)+v(x,t)\cdot\grad\zeta(x,t)&=0,\qquad x\in\R^2,\, t\in[0,\infty),\\
 \zeta(x,0)&=\zeta_0(x),
\end{split}
\end{align}
where $\zeta$ is an active scalar driven by the incompressible velocity $v(x,t)$:
\begin{equation}\label{incompre}
\nabla\cdot v(x,t)=0.
\end{equation}
Depending upon our choice of the relation between the velocity and the scalar we will obtain from this system both the Surface Quasi-geostrophic equation (SQG) for sharp fronts and the Muskat problem. In this paper we will present a unified method to establish the absence of splash singularities for both of these systems in different scenarios. Specifically we will show that the dynamics of a smooth contour can not cause an intersection at a single point.

We obtain the SQG equations from the system \eqref{ITE} and  \eqref{incompre} by expressing the velocity $v$ in terms of a stream function
$$
v = \nabla^\perp \psi = (-\partial_{x_2}\psi,  \partial_{x_1}\psi ),
$$
where the function $\psi$ satisfies $\zeta = - (-\Delta)^{1/2} \psi$. Here $(-\Delta)^{1/2}$ is the Zygmund operator defined on the Fourier side by $$\widehat{(-\Delta)^{1/2}}=|\xi|.$$ This can be shown to be equivalent to the condition
\begin{equation}\label{SQG}
v(x,t)=(-R_2\zeta(x,t),R_1\zeta(x,t)),
\end{equation}
which relates the temperature to the velocity by means of the Riesz transforms $R_1$ and $R_2$.

The SQG system is physically important as a model of atmospheric turbulence and oceanic flows (see, e.g., \cite{CMT}, \cite{Held}, \cite{CNS} and the references therein). This equation is derived in the situation of small Rossby and Ekman numbers and constant potential vorticity \cite{Gill}, where the scalar $\zeta$ is the evolution over time of the temperature of the fluid. SQG has been the subject of a large number of studies from different points of view.   Underlying its mathematical interest are its  strong analogies with the three-dimensional Euler equations (see \cite{CMT} and \cite{Majda} for these discussions).   A very actively studied question for this system has been the formation of singularities in finite time for smooth initial data (see, e.g., \cite{Ohkitani}, \cite{D}, \cite{DC2}, \cite{Hou}, \cite{CCW}, \cite{CLSTW} and references therein).

The SQG system has been furthermore used as a mathematical model in the meteorological process of frontogenesis. Here the dynamics of hot and cold fluids are studied in the context of the formation and time evolution of weather sharp fronts where the temperature exhibits discontinuity jumps (further information can be found in \cite{CMT} and the references therein). In light of this interest, Rodrigo \cite{Rodrigo} studied the case in which the initial temperature takes two different constant values on complementary domains
\begin{equation}\label{initialdata}
\zeta_0(x)=\left\{\begin{array}{cl}
                    \zeta^1,& x\in\Omega_0,\\
                    \zeta^2,& x\in\R^2\smallsetminus\Omega_0,
                 \end{array}\right.
\end{equation}
where $\zeta^1 \neq \zeta^2$. The initial data represent sharp fronts and the interest is in their dynamics which evolve by SQG. The transport character of the equation \eqref{ITE} shows that the temperature as it evolves in time should have the form
\begin{equation}\label{chi}
\zeta(x,t)=\left\{\begin{array}{cl}
                    \zeta^1,& x\in\Omega(t),\\
                    \zeta^2,& x\in\R^2\smallsetminus\Omega(t).
                 \end{array}\right.
\end{equation}
In this formulation $\Omega(t)$ is a moving domain.
Then a contour dynamics problem is obtained by considering the time evolution of the free boundary $\partial\Omega(t)$.
For the equations \eqref{ITE}, \eqref{incompre} and \eqref{SQG} the SQG equation for the evolution of a sharp front is then given by:
\begin{equation}\label{Ifront}
x_t(\al)=\frac{\zeta^2-\zeta^1}{2\pi}\int_{\T}\frac{x'(\al)-x'(\beta)}{|x(\al)-x(\beta)|}d\beta.
\end{equation}
Here the boundary is parameterized by the closed one-to-one curve $x(\al, t)$:
\begin{equation*}
\partial\Omega(t)=\{x(\al,t)=(x_1(\al,t),x_2(\al,t))\,:\,\al\in[-\pi,\pi]=\T\},
\end{equation*}
which satisfies the chord-arc condition (see \cite{G} for a detailed derivation of the contour equation in this form). Above the subscript $t$ and the prime notation denote the partial derivatives in time and $\alpha$ respectively. In \eqref{Ifront} the time dependence is disregarded for notational simplicity.

Then fundamental questions to study are the existence of front type solutions and the possible singularity formation in the evolution of $\partial\Omega(t)$. These issues are comparable to the vortex-patch problem for the 2-D Euler equations (see \cite{bertozzi-Constantin} and \cite{Majda}), but the SQG front system is more singular (see \cite{Rodrigo} for more details on this discussion).

Local-in-time existence and uniqueness in this situation was proven in \cite{Rodrigo} for $C^\infty$ contours using Nash-Moser arguments.  This tool was used because the operator involved in the contour equation is considerably singular, it loses more than one derivative (see equation \eqref{Ifront} above).  In \cite{G} the result was extended within the chain of Sobolev spaces because of several cancelations.  See also \cite{FR} for a proof of local existence for analytic contours. In \cite{CFMR} numerical simulations indicate the possibility of singularity formation on the free boundary.  More specifically, initial data was shown  in which the curvature blows up numerically due to two branches of the fluid interphase collapsing in a single point in a self-similar way. Therefore this work provides an interesting stable scenario for a possible singularity formation. Recently there has been an active interest in the study of almost sharp front type weak solutions of the SQG equation (see \cite{FR2} and the references therein for more details).

We next discuss the Muskat problem; this system models the physical scenario of multiple fluids with different characteristics in porous media. Specifically we will study the dynamics of interphases between fluids which are immiscible and incompressible. To derive the equations of the Muskat problem, the system \eqref{ITE}-\eqref{incompre} is used and we choose the velocity $v(x,t)$ to satisfy
 Darcy's law:
\begin{equation}\label{IPM}
v(x,t)=-\nabla p(x,t)-(0,\zeta(x,t)).
\end{equation}
Above the scalar $p(x,t)$ is the pressure and in this situation $\zeta$ is the scalar density.  Also the acceleration due to gravity and the viscosity of the fluid are set to unity to simplify the notation. Then the system \eqref{ITE}-\eqref{incompre} turns out to be the conservation of mass, which together with \eqref{IPM} yields the incompressible porous media equation (IPM): \cite{Otto}, \cite{CGO}, \cite{Lazlo}.  By considering a solution of the form \eqref{chi}, then the interphase $\partial\Omega(t)$ is a free boundary and it describes the density jump between each fluid. The evolution equation is given by
\begin{equation}\label{IMuskat}
x_t(\al)=\frac{\zeta^2-\zeta^1}{2\pi}\int\frac{(x_1(\al)-x_1(\beta))(x'(\al)-x'(\beta))}{|x(\al)-x(\beta)|^2}d\beta.
\end{equation}
(see \cite{ADCPM} for the whole derivation). Above $\alpha,\beta\in\T$ for closed contours and  $\alpha,\beta\in\R$ for the asymptotically flat case when
$x(\alpha,t)-(\alpha,0)\to 0$ as $\alpha\to \infty$.  Further $\alpha,\beta\in\R$ for periodic curves in the $x_1$ direction when
$x(\alpha+(2\pi,0),t)=x(\alpha,t)$. The integral in \eqref{IMuskat} is understood as a principal value when that is necessary.

The Muskat problem is a classical well-established problem \cite{Muskat}. It has been highly studied  in particular because of strong similarities with the contour dynamics of fluids in Hele-Shaw cells \cite{ST}. For both of these completely different physical scenarios
it is therefore possible to reach similar conclusions.  These problems can be studied in the situations of fluids with different viscosities \cite{SCH} and with surface tension
effects \cite{Esch2}. Notice that the formulation \eqref{IMuskat} above describes the case in which the viscosities and pressures are equal across the interphase (and \eqref{IMuskat} is in the situation with no surface tension).

The Muskat problem has been shown to exhibit instabilities and ill-posedness in several situations (see for instance \cite{Otto}, \cite{SCH}, \cite{Esch2}, \cite{Lazlo}). For the situation that we are studying in this paper, e.g. the contour evolution system \eqref{IMuskat}, the instabilities in the system will appear when the heavy fluid lies on top of the light one (see \cite{DY}).

When the light fluid lies above the heavy fluid this is called the stable scenario and in this case the system has been shown to be well-posed \cite{DY}.  More generally, for the Muskat problem, the well-posedness condition amounts to the positivity of the difference of the gradient of the pressure jump at the interphase in Darcy's law \eqref{IPM} (see \cite{Ambrose}).  This condition must hold for the initial data in order for the system to be well-posed  \cite{ADY}.  It is known in the literature as the Rayleigh-Taylor sign condition \cite{ST}, \cite{Lannes}. The stable framework gives rise to global-existence results for  initial data with small norms \cite{SCH}, \cite{Esch2}, \cite{ccgs}, \cite{Beck}. On the other hand, global-existence can be false for certain scenarios with large initial data.  In \cite{ADCPM} it was proven that there exist initial data in the stable regime for equation \eqref{IMuskat} such that the solution turns to the unstable regime in finite time. This interface was initially a smooth stable graph (with the heavier fluid below), but later it enters into an unstable regime.  In other words the interphase is transformed into a non-graph in finite time: when this happens we say that the interphase ``turns-over''.  The particular significance of a turn-over is that the Rayleigh-Taylor condition breaks down.  At some branch in the interphase it is possible to localize the heavy fluid on top of the lighter one.  Then at a later time the regularity of the contour breaks down \cite{ADCP}, i.e., the Muskat problem develops a singularity in finite time starting from regular stable initial data.

We briefly discuss the 2D water waves problem, which is another incompressible fluid interphase dynamics equation.  This system can be given by \eqref{ITE}-\eqref{incompre} together with the 2D density variable Euler equations:
\begin{equation}\label{Euler}
\zeta(x,t)(v_t+v\cdot\grad v)(x,t)=-\grad p(x,t)-(0,\zeta(x,t)).
\end{equation}
We have solutions to this system in the form of \eqref{chi}, which establishes the evolution of a free boundary given by air, $\zeta^1=0$, and
water, $\zeta^2=1$, governed by the gravity force. The velocity is assumed to be rotationally free on each side of $\partial\Omega(t)$ but concentrated on the moving interphase as a delta distribution: \begin{equation}\label{wdelta}
\nabla^\perp\cdot v(x,t)=\omega(\al,t)\delta(x=x(\al,t)).\end{equation}
 There is a large mathematical literature  on the 2D water waves problem (see \cite{Lannes} and the references therein).
This system has been shown to be well-posed if the Rayleigh-Taylor condition is satisfied initially \cite{Wu}.  There are recent global-in-time results for small initial data (see \cite{Wu2}, \cite{Ionescu}, \cite{Alazard} and the references therein).  On the other hand, for large initial data with an ``overturning shape'' the system develops finite time splash singularities \cite{ADCPJ}, \cite{ADCPJ2}, \cite{CS}.    More precisely there exists a family of
initial data satisfying the chord-arc condition such that the interface $x(\al,t)$ from the solution of the system \eqref{ITE}, \eqref{incompre}, \eqref{chi}, \eqref{Euler} and \eqref{wdelta} satisfying  $\zeta^1=0$ touches itself at a single point at time $t_s>0$ meanwhile $x(\al,t_s)$ is smooth. In particular, the curvature is finite.  We would also like to mention recent developments by C. Fefferman, A. Ionescu and V. Lie on the absence of splash singularities for two incompressible fluids.

With the results below we prove that in order to have a pointwise collapse, the second derivative, and therefore the curvature, has to blow-up. Splash singularities turn out to be false for the SQG sharp fronts and the Muskat problem. This phenomena was observed numerically in \cite{CFMR}, where computer solutions of the SQG sharp front system exhibit pointwise collapse and the curvature blows-up at the same finite time.

We also improve the result in \cite{DP3}, where it is shown that a positive volume of fluid between the contours can not be ejected in finite time. That result is proved by showing that the velocity is bounded \cite{DC2} for the Muskat problem for smooth contours. The velocity can be related to the density using singular integral operators with even kernels \eqref{IPM}, Then the fact that $\zeta$ is given by a step function \eqref{chi} allows one to show that $v$ is in $L^\infty$.  There is a cancellation that is used to establish that $v$ is bounded which was previously obtained by Bertozzi-Constantin in \cite{bertozzi-Constantin}.  They applied it to the 2D vortex-patch problem to prove global regularity.  The present work contributes the information that the level set can not collapse even pointwise.

The pointwise collapse of smooth level sets, and therefore splash singularities, for regular solutions of SQG, IPM and general active scalar equations have been extensively studied (see e.g. \cite{DC2}, \cite{CCW}).  Although, for initial data that is not necessarily a sharp front \eqref{initialdata}, the situation could be a priori less singular and the problem is still open.

We will explain the proof of our results first for the multi-phase Muskat problem.  Our reasoning is two fold.  First the Muskat scenario that we present is well-posed \cite{DP3} and there is no Rayleigh-Taylor instabilities \cite{Esch3}. Second, the proof in this case will appear more clearly.  We will consider fluids which have three different constant values for the density:
\begin{equation*}
 \zeta(x_1,x_2,t)=\left\{\begin{array}{cl}
                    \zeta^1\quad\mbox{in}&\{x_2>f(x_1,t)\},\\
                    \zeta^2\quad\mbox{in}&\{f(x_1,t)>x_2>g(x_1,t)\},\\
                    \zeta^3\quad\mbox{in}&\{g(x_1,t)>x_2\},
                 \end{array}\right.
\end{equation*}
where we suppose that $f(x_1,t)>g(x_1,t)$, and that the two dynamic surfaces, which are defined by
$x_2= f(x_1,t)$ and $x_2= g(x_1,t)$, can be parametrized as a graph at time $t=0$. The constant densities satisfy
\begin{equation}\label{den}
\zeta^1<\zeta^2<\zeta^3.
\end{equation}
This keeps us in the stable situation. Furthermore we work in the situation where
$$
\lim_{x_1\to\infty}f(x_1,t)=f_\infty> g_\infty=\lim_{x_1\to\infty}g(x_1,t).
$$
 Then our result can be stated as follows:

 \begin{thm}\label{MuskatTHM}
 Suppose the free boundaries $f(\al,t)$ and $g(\al,t)$ are smooth for $\alpha\in \R$ and $t\in [0,T)$  with $T>0$ arbitrary.  Define the distance:
\begin{equation}\label{minMUSKAT}
0<S(t)=\min_{\al\in \R}(f(\al,t)-g(\al,t))\ll \min\{f_\infty-g_\infty,1\}.
\end{equation}
Then the following uniform lower bound for $t\in [0,T)$ holds:
\begin{equation}\label{MUSKATlower}
S(t)\geq \exp\left(\ln(S(0)) \exp\left({\int_0^t C(f,g)(s)ds}\right)\right).
\end{equation}
Here $C(f,g)$ is a smooth function of $\|f''\|_{L^\infty}+\|g''\|_{L^\infty}$ and $\|f\|_{L^\infty}+\|g\|_{L^\infty}$ which is defined in \eqref{cMUSKATdef} below.
\end{thm}

After proving Theorem \ref{MuskatTHM} we will extend these results to the SQG sharp front system based on the previous approach which was used for the Muskat problem.

\begin{thm}
Consider a smooth curve $x(\alpha, t)$ which is a solution to the sharp front SQG system for $t\in [0,T)$ with $T>0$ arbitrary.  Let $S(t)>0$ be defined as the minimum distance between two different branches of the interphase which are approaching each other as $t\to T^+$. Then $S(t)$ is bounded below by an explicitly computable positive function which goes to zero double exponentially fast for $t$ traveling to infinity.
\end{thm}

Finally at the end of this paper we will show additional scenarios where our result will hold such as the multi-phase SQG system. For Muskat we also consider the cases of closed contours and overturning shaped interphases. Although in those situations Rayleigh-Taylor instabilities appear and the interphases have to be analytic in order for there to be bona fide solutions (see \cite{ADCP}).

\section{The Multi-phase Muskat Problem}

The contour equation for the multiphase Muskat problem can be written as
\begin{align}
\begin{split}\label{ec2i1}
\D f_t(\al)&=\int_{\R}(\zeta^{21}K(f,f)+\zeta^{32}K(f,g))(\alpha,\beta)d\beta,\\
\D g_t(\al)&=\int_{\R}(\zeta^{32}K(g,g)+\zeta^{21}K(g,f))(\alpha,\beta)d\beta,
\end{split}
\end{align}
where $\zeta^{21} \eqdef (\zeta^2-\zeta^1)/(2\pi)$, $\zeta^{32} \eqdef (\zeta^3-\zeta^2)/(2\pi)$,
$$
K(f,g)(\alpha,\beta)=\frac{\beta\delta_\beta(f',g')(\alpha)}{\beta^2+(\delta_\beta(f,g)(\alpha))^2},
$$
$$\delta_\beta(f,g)(\alpha)=f(\al)-g(\al-\beta),$$ and for simplicity we denote
$$
\delta_\beta f(\alpha)=f(\al)-f(\al-\beta).
$$
We remark that it is possible to recover \eqref{IMuskat} by taking $\zeta^{32}=g(\al)=0$ and $x(\al,t)=(\al,f(\al,t))$ (see \cite{DP3} for a detailed derivation of this equation).

We next check the evolution of \eqref{minMUSKAT} and denote $\al_t\in\R$ such that $S(t)=f(\al_t,t)-g(\al_t,t)$. We use the Rademacher theorem to obtain that $S(t)$ is differentiable almost everywhere and that
$S_t(t)=f_t(\al_t,t)-g_t(\al_t,t)$ (see \cite{CE} and \cite{AD} for the whole argument). We plug this identity into \eqref{ec2i1} to split the integration regions as
\begin{align}\notag
\begin{split}
S_t(t)
&=\int_{|\beta|<S(t)}d\beta+\int_{S(t)<|\beta|<1}d\beta+\int_{|\beta|>1}d\beta  \\
&=I+II+III.
\end{split}
\end{align}
For the first integral we bound the kernels $K$ in absolute value using the crucial identity
$$
f^\prime(\al_t,t)=g^\prime(\al_t,t),
$$
to find
$$
I\leq 2(\zeta^{21}\|f''\|_{L^\infty}+\zeta^{32}\|g''\|_{L^\infty})\int_{|\beta|<S(t)}1 d\beta,
$$
and therefore we obtain
$$
I\leq C(\|f''\|_{L^\infty}+\|g''\|_{L^\infty})S(t).
$$
For the second integral we split further $II=\zeta^{21}II_1+\zeta^{32}II_2$. We will show how to deal with $II_1$ and observe that $II_2$ is analogous. Notice that we have
\begin{equation}
\begin{split}
II_1=&\int_{S(t)<|\beta|<1}\!\!\!\!\!\!\!\!\!\!\!\!\!\!\!\!\!\frac{\beta\delta_\beta f'(\al_t) [ (\delta_\beta(g,f)(\alpha_t))^2-(\delta_\beta f(\alpha_t))^2]}{D(g,f,\beta)}d\beta,
\end{split}
\end{equation}
where the denominator is given by
$$
D(g,f,\beta) \eqdef [\beta^2+(\delta_\beta f(\alpha_t))^2][\beta^2+(\delta_\beta(g,f)(\alpha_t))^2].
$$
We split further
\begin{equation}\label{bob}
\begin{split}
II_1=&-\int_{S(t)<|\beta|<1}\frac{\beta \delta_\beta f'(\al_t) S(t)
\delta_\beta(g,f)(\alpha_t) }{D(g,f,\beta)}d\beta
\\
&-\int_{S(t)<|\beta|<1}\frac{\beta \delta_\beta f'(\al_t) S(t)
\delta_\beta f(\alpha_t) }{D(g,f,\beta)}d\beta,
\end{split}
\end{equation}
to obtain
$$
II_1\leq 2\|f''\|_{L^\infty}S(t)\int_{S(t)<|\beta|<1}|\beta|^{-1}d\beta.
$$
The last calculation yields
$$
II_1\leq -2\|f''\|_{L^\infty}S(t)\ln S(t).
$$
The term $II_2$ can be estimated similarly and we obtain
$$
II\leq -C(\|f''\|_{L^\infty}+\|g''\|_{L^\infty})S(t)\ln S(t).
$$
For the last term $III$ we arrange the terms as in \eqref{bob} to find
\begin{equation}\notag
\begin{split}
III\leq C(f,g) S(t)\int_{|\beta|>1}|\beta|^{-2}d\beta,
\end{split}
\end{equation}
with
\begin{equation}\label{cMUSKATdef}
C(f,g)=C(\|f''\|_{L^\infty}+\|g''\|_{L^\infty})(\|f\|_{L^\infty}+\|g\|_{L^\infty}+1).
\end{equation}
Collecting all of the previous estimates we obtain that
$$
S_t(t)\geq C(f,g)S(t)\ln S(t).
$$
A further time integration yields \eqref{MUSKATlower}.
Notice that $\ln(S(0)) < 0$.  Thus $S(t)$ can not go to zero in finite time.

\section{SQG sharp front}

For the SQG sharp front equation we choose the parametrization for the contour equation that yields the equation
\begin{equation}\label{front}
x_t(\al)=\int_{\T}\frac{\delta_\beta x'(\al)}{|\delta_\beta x(\al)|}d\beta,
\end{equation} where we take $\zeta^2-\zeta^1=2\pi$ for the sake of simplicity.
We now assume without loss of generality that the pointwise approaching ``splash'' is going to take place in a small ball $B$ of radius $\epsilon_0/2$ and center $(0,0)$. The two branches of the interfaces will be approaching horizontally so that they are represented by $(\al,f(\al,t))$ and $(\al,g(\al,t))$ inside $2B$ with $f>g$. We then find for the chart $x(\al)=(\al,f(\al))$ for $\al\in(-\epsilon_0,\epsilon_0)$ the equation
\begin{equation}\notag
\begin{split}
f_t(\al)=&\int_{-\epsilon_0}^{\epsilon_0}\frac{\delta_\beta f'(\al)}{\sqrt{\beta^2+(\delta_\beta f(\al))^2}}d\beta  \\
&+\int_{\epsilon_0}^{-\epsilon_0}\frac{\delta_\beta (f',g')(\al)}{\sqrt{\beta^2+(\delta_\beta(f, g)(\al))^2}}d\beta+R(f),
\end{split}
\end{equation}
where $R(f)$ is the remainder in the integral equation in the second component in \eqref{front}.
For the chart $x(\al)=(\al,g(\al))$ we similarly have
\begin{equation}\notag
\begin{split}
g_t(\al)=&\int_{-\epsilon_0}^{\epsilon_0}\frac{\delta_\beta (g',f')(\al)}{\sqrt{\beta^2+(\delta_\beta (g,f)(\al))^2}}d\beta  \\
& +\int_{\epsilon_0}^{-\epsilon_0}\frac{\delta_\beta g'(\al)}{\sqrt{\beta^2+(\delta_\beta g(\al))^2}}d\beta+R(g),
\end{split}
\end{equation}
and $R(g)$ is again the remainder given by \eqref{front}. We now define as before
$$
S(t) \eqdef \min_{[-\epsilon_0,\epsilon_0]}(f(\al,t)-g(\al,t))=f(\al_t,t)-g(\al_t,t),
$$
with $\al_t\in (-\epsilon_0/2,\epsilon_0/2)$. We proceed as in the previous section to follow $S_t(t)$ for almost every $t$. We find that the integrals above can be handled in a similar manner except for $R(f)$ and $R(g)$. For these remainder terms we can choose $\epsilon_0$ small enough so that $R(f)=O(S(t))$ and $R(g)=O(S(t))$. In fact
$$
|R(f)|\leq \frac{\|x''\|_{L^\infty}}{c_{CA}}\int_{\T\smallsetminus[-\epsilon_0,\epsilon_0]}d\beta,
$$
where $c_{CA}>0$ is the chord-arc constant of the curve outside the ball $B$:
\begin{equation*}
|\delta_{\beta}x(\al)|\geq c_{CA} |\beta|,\,\,\al\in[-\epsilon_0/2,\epsilon_0/2],\,\,\beta\in\T\smallsetminus[-\epsilon_0,\epsilon_0].
\end{equation*}
The analogous estimate  for $R(g)$ follows similarly.  We can thus obtain
$$
S_t(t)\geq C(x)S(t)\ln S(t),
$$
where $C(x)=C(\|x''\|_{L^\infty},c_{CA},\epsilon_0)$. We therefore again control the size of $S(t)$ from below by double exponential time decay.

\section{Additional scenarios for Muskat and SQG}
This analysis also works for the the multi-phase SQG sharp front system. In that case the equations for the $2\pi$-periodic contours $f(\al,t)$ and $g(\al,t)$ are given by
\begin{align}
\begin{split}\notag 
\D f_t(\al)&=\int_{\T}(\zeta^{21}\Sigma(f,f)+\zeta^{32}\Sigma(f,g))(\alpha,\beta)d\beta,
\\
\D g_t(\al)&=\int_{\T}(\zeta^{32}\Sigma(g,g)+\zeta^{21}\Sigma(g,f))(\alpha,\beta)d\beta.
\end{split}
\end{align}
Above $\zeta^{21}=(\zeta^2-\zeta^1)/(2\pi)$ and $\zeta^{32}=(\zeta^3-\zeta^2)/(2\pi)$ with no order needed in the size of $\zeta^1$, $\zeta^2$ and $\zeta^3$ as in \eqref{den} because there are no instabilities for SQG. The kernel $\Sigma(f,g)(\al,\beta)$ behaves like
$$
\Sigma(f,g)(\al,\beta)=\frac{\delta_{\beta}(f',g')(\al)}{\sqrt{\beta^2+(\delta_{\beta}(f,g)(\al))^2}},
$$
for $\beta$ close to $0$ and $f(\al)$ close to $g(\al-\beta)$. Hence the same  approach as  described previously for SQG follows.

We end by proposing two additional scenarios. These are closed and overturning shaped contours for the Muskat equation \eqref{IMuskat}. In those cases the same results can be shown as for the SQG sharp fronts. But because of the Rayleigh-Taylor instabilities the solutions to the interphase equations have to be analytic in order to make rigorous mathematical sense \cite{ADCPM}.

\subsection*{{\bf Acknowledgements}}
\smallskip
FG was partially supported by MCINN grant MTM2011-26696 and Ram\'on y Cajal program (Spain). RMS was partially supported by the NSF grant DMS-1200747 and an Alfred P. Sloan Foundation Research Fellowship.

%



\bigskip


\begin{quote}
\begin{tabular}{ll}
\textbf{Francisco Gancedo}\\
{\small Departamento de An\'{a}lisis Matem\'{a}tico}\\
{\small Universidad de Sevilla}\\
{\small Tarfia s/n, 28006 Sevilla, Spain}\\
{\small Email: fgancedo@us.es}
\end{tabular}
\end{quote}

\begin{quote}
\begin{tabular}{ll}
\textbf{Robert M. Strain}\\
{\small Department of Mathematics}\\
{\small University of Pennsylvania}\\
{\small David Rittenhouse Lab}\\
{\small 209 South 33rd Street, Philadelphia, PA 19104, USA} \\ 
{\small Email: strain@math.upenn.edu}
\end{tabular}
\end{quote}


\begin{thebibliography}{99}

\bibitem{Alazard} T. Alazard and J.M. Delort. Global solutions and asymptotic behavior for two dimensional gravity water waves. ArXiv:1305.4090 (2013).

\bibitem{Ambrose} D.M. Ambrose. Well-posedness of two-phase Hele-Shaw flow without surface tension. European J. Appl. Math. 15, no. 5, 597-607 (2004).

\bibitem{Beck} T. Beck, P. Sosoe and P. Wong. Duchon-Robert solutions for the Rayleigh-Taylor and Muskat problems. ArXiv:1209.1113v3 (2013).

\bibitem{bertozzi-Constantin} A.~L. Bertozzi and P. Constantin. Global regularity for vortex patches.
\emph{Comm. Math. Phys.} 152 (1): 19--28, 1993.

\bibitem{ADCP} A. Castro, D. C\'ordoba, C. Fefferman and F. Gancedo. Breakdown of smoothness for the Muskat problem.
Arch. Ration. Mech. Anal., 208, no. 3, 805-909 (2013).

\bibitem{ADCPJ} A. Castro, D. C\'ordoba, C. Fefferman, F. Gancedo and J. G{\'o}mez-Serrano. Splash singularity for water waves.
Proc. Natl. Acad. Sci., 109, no. 3, 733-738 (2012).

\bibitem{ADCPJ2} A. Castro, D. C\'ordoba, C. Fefferman, F. Gancedo and J. G{\'o}mez-Serrano. Finite time singularities for the free boundary incompressible Euler equations. \emph{To appear in Annals of Math} (2013).

\bibitem{ADCPM} A. Castro, D. C\'ordoba, C. Fefferman, F. Gancedo and M. L\'opez-Fern\'andez. Rayleigh-Taylor breakdown for the Muskat
problem with applications to water waves. Annals of Math, 175, no. 2, 909-948 (2012).

\bibitem{CCW} D. Chae, P. Constantin and J. Wu. Deformation and symmetry in the inviscid SQG and the 3D Euler equations. J. Nonlinear Sci. 22 no. 5, 665-688, (2012).

\bibitem{CE} A. Constantin and J. Escher. Wave breaking for nonlinear nonlocal shallow water equations.
\emph{Acta Math.}, 181, 229-243 (1998).

\bibitem{ccgs} P. Constantin, D. C\'ordoba, F. Gancedo and R.M. Strain. On the global existence  for
the Muskat problem. \emph{J. Eur. Math. Soc.}, 15, 201-227 (2013).

\bibitem{CLSTW} P. Constantin, M.C. Lai, R. Sharma, Y.H. Tseng, J. Wu. New numerical results for the surface quasigeostrophic equation.
\emph{J. Sci. Comput.} 50, 1-28 (2012).

\bibitem{CMT} P.~Constantin, A.~J. Majda, and E.~Tabak. \newblock Formation
of strong fronts in the 2-{D} quasigeostrophic thermal active scalar.
\newblock \emph{Nonlinearity}, 7:1495--1533, 1994.

\bibitem{CNS} P. Constantin, Q. Nie, and N. Schorghofer. Nonsingular surface quasi-geostrophic flow.
Phys. Lett. A, 241, no. 3, 168-172, 1998.

\bibitem{D} D. C\'ordoba. Nonexistence of simple hyperbolic blow-up for the quasi-geostrophic equation. Ann. of Math. (2) 148, no. 3, 1135-1152, (1998).

\bibitem{AD} A. C\'ordoba, D. C\'ordoba. A pointwise estimate for fractionary derivatives with applications to partial differential equations. \emph{Proc. Natl. Acad. Sci. USA} 100, 26,
(2003), 15316-15317.

\bibitem{ADY} A. C\'ordoba, D. C\'ordoba and F. Gancedo. Interface evolution: the Hele-Shaw and Muskat problems.
 \emph{Annals of Math.}, 173, 1, (2011), 477-542.

\bibitem{DC2} D. C\'ordoba, Ch. Fefferman. Scalars convected by a two-dimensional incompressible flow. Comm. Pure Appl. Math. 55, no. 2, 255-260 (2002).

\bibitem{CFMR} D. C\'ordoba, M. A. Fontelos, A. M. Mancho and J. L. Rodrigo. Evidence of singularities
for a family of contour dynamics equations. \newblock\emph{Proc. Natl. Acad. Sci. USA} 102,
5949--5952, 2005.

\bibitem{DY} D. C\'ordoba and F. Gancedo. Contour dynamics of incompressible 3-D fluids
in a porous medium with different densities. \emph{Comm. Math.
Phys.} 273, 2, (2007), 445-471.

\bibitem{DP3} D. C\'ordoba and F. Gancedo. Absence of squirt singularities for the multi-phase Muskat problem.
\emph{Comm. Math. Phys.}, 299, 2, (2010), 561-575.

\bibitem{CGO} D. C\'ordoba, F. Gancedo and R. Orive. Analytical behavior of 2D incompressible flow in porous media. J. Math. Phys., 48, no. 6 (2007).

\bibitem{CS} D. Coutand and S. Shkoller. On the finite-time splash and splat singularities for the 3-D free-surface Euler equations.
\emph{to appear in Commun. Math. Phys.} (2013).

\bibitem{Hou} J. Deng, T.Y. Hou, R. Li, X. Yu. Level set dynamics and the non-blowup of the 2D quasi-geostrophic
equation. Methods Appl. Anal. 13, 157-180 (2006).

\bibitem{Esch2} J. Escher and B.V. Matioc. On the parabolicity of the Muskat problem: Well-posedness, fingering,
and stability results. Z. Anal. Anwend. 30, no. 2, 193-218 (2011).

\bibitem{Esch3} J.  Escher, A.V. Matioc and B.V. Matioc. A generalized Rayleigh-Taylor condition for the Muskat problem. Nonlinearity 25, no. 1, 73-92 (2012).

\bibitem{FR} C. Fefferman, J.L. Rodrigo. Analytic sharp fronts for the surface quasi-geostrophic equation. Comm. Math. Phys. 303, no. 1, 261-288 (2011).

\bibitem{FR2} C. Fefferman, J.L. Rodrigo. Almost sharp fronts for SQG: the limit equations. Comm. Math. Phys. 313, no. 1, 131-153 (2012).

\bibitem{G} F. Gancedo. Existence for the $\alpha$-patch model and the QG sharp front in Sobolev spaces. Adv.Math. 217(6),
2569-2598, (2008).

\bibitem{Gill} A.E. Gill. Atmosphere-Ocean Dynamics. Academic Press, New York (1982).

\bibitem{Held} I. Held, R. Pierrehumbert, S. Garner and K. Swanson. Surface quasi-geostrophic dynamics. J. Fluid Mech.
282, 1-20 (1995).

\bibitem{Ionescu} A.D. Ionescu and F. Pusateri. Global solutions for the gravity water waves system in 2d. ArXiv:1303.5357 (2013).

\bibitem{Lannes} D. Lannes. The Water Waves Problem: Mathematical Analysis and Asymptotics. American Mathematical Society (2013).

\bibitem{Majda} A. Majda and A. Bertozzi. Vorticity and Incompressible Flow. Cambridge University Press, Cambridge
(2002).

\bibitem{Muskat} M. Muskat. Two Fluid systems in porous media. The encroachment of water into an oil
sand. Physics, 5, 250-264 (1934).

\bibitem{Ohkitani} K. Ohkitani, M. Yamada. Inviscid and inviscid-limit behavior of a surface quasigeostrophic Flow.
Phys. Fluids 9, 876-882 (1997).

\bibitem{Otto} F. Otto. Evolution of microstructure in unstable porous media flow: a relaxational approach.
\emph{Commun. Pure Appl. Maths} 52, 873-915 (1999).

\bibitem{Rodrigo} J.L. Rodrigo. On the Evolution of Sharp Fronts for the
Quasi-Geostrophic Equation. \emph{Comm. Pure and Appl. Math.}, 58: 0821-0866, 2005.

\bibitem{ST} P.G. Saffman and G. Taylor. The penetration of a fluid into a porous medium or Hele-Shaw cell
containing a more viscous liquid. Proc. R. Soc. London, Ser. A 245, 312-329, 1958.

\bibitem{Lazlo} L. Jr. Sz\'ekelyhidi. Relaxation of the incompressible porous media equation. Ann. Sci. \'Ec. Norm. Sup\'er. (4) 45, no. 3, 491-509, 2012.

\bibitem{SCH} M. Siegel, R. Caflisch and S. Howison. Global
Existence, Singular Solutions, and Ill-Posedness for the Muskat
Problem. \emph{Comm. Pure and Appl. Math.}, 57, (2004), 1374-1411.

\bibitem{Wu} S. Wu. Well-posedness in Sobolev spaces of the full water wave problem in 2-D. Invent.
Math., 130(1):39-72, 1997.

\bibitem{Wu2} S. Wu. Almost global wellposedness of the 2-D full water wave problem. Invent. Math.,
177(1):45-135, 2009.

\end{thebibliography}
\end{document}